\documentclass[%
 aip,
 jmp,%
 amsmath,amssymb,
 reprint,%
]{revtex4-1}

\usepackage{mathtools}
\usepackage{floatrow}
\usepackage{graphicx}
\usepackage{dcolumn}
\usepackage{bm}
\usepackage[inline]{enumitem}

\def\R{\mathbb{R}}
\def\Z{\mathbb{Z}}

\def\mc{\mathcal}

\numberwithin{equation}{section}
\DeclareMathOperator{\sinc}{sinc}
\newtheorem{conjecture}{Conjecture}

\begin{document}

\title[]{Jumping champions and prime gaps using information-theoretic tools}
\author{Nicholas Pun}
\affiliation{Department of Applied Mathematics, University of Waterloo}
\author{Robert T.W. Martin}
\affiliation{Department of Mathematics and Applied Mathematics, University of Cape Town}
\author{Achim Kempf}
\affiliation{Departments of Applied Mathematics and Physics, University of Waterloo}
\affiliation{Institute for Quantum Computing, University of Waterloo}
\date{\today}

\begin{abstract}
We study the spacing of the primes using methods from information theory. In information theory, the equivalence of continuous and discrete representations of information is established by Shannon sampling theory. Here, we use Shannon sampling methods to construct continuous functions whose varying bandwidth follows the distribution of the prime numbers. The Fourier transforms of these signals spike at frequently occurring spacings between the primes. We find prominent spikes, in particular, at the primorials. Previously, the primorials have been conjectured to be the most frequent gaps between subsequent primes, the so-called ``jumping champions". Here, we find a foreshadowing of the primorial's role as jumping champions in the sense that Fourier spikes for the primorials arise much earlier on the number axis than where the primorials in question are expected to reign as jumping champions. 
\end{abstract}

\maketitle

\section{Introduction}
The gaps between the primes possess intriguing structural properties and have led to a number of important as yet unproven conjectures, such as the Hardy-Littlewood $k$-tuple conjecture\cite{hardy1923} of 1923. More recently, in 1999, Odlyzko, Rubinstein and Wolf \cite{jumpingchampions} published a conjecture concerning so-called jumping champions. For any $t> 0$, the jumping champion is defined as the integer, $g$, that is the most frequently occurring gap between any two successive primes less than or equal to $t$. 
The jumping champions conjecture states:
\begin{conjecture} \label{conj:jc}
The jumping champions are 4 and the primorials, i.e.,  2, 6, 30, 210, ....
\end{conjecture}
Here, the $n$'th \textit{primorial} is the product of the first $n$ primes. In 2012, Goldston and Ledoan proved that a version of the Hardy-Littlewood $k$-tuple conjeture for prime pairs and triples implies that all sufficiently large jumping champions are primorials, and any sufficiently large primorials are jumping champions over a long range of $t>0$, see \cite{jumpingchampions2}. In particular, they provide estimates on ranges of $t$ for which a given primorial is the jumping champion for $[0,t]$. For example, the primorial 210 is expected to reign as jumping champion in (roughly) the interval $[10^{487},10^{2607}]$. See \cite{Wolf1999,Wolf2011,Ares2006,Szpiro2004} for related results and investigations. 

The magnitude of these numbers would appear to preclude numerical studies. 
However, as we will show, intriguing evidence for the importance of the primorials in the distribution of prime gaps can be obtained numerically through the use of the information-theoretic tools of Shannon sampling theory. In information theory, Shannon sampling constructively establishes the equivalence between discrete and continuous representations of information\cite{shannon,Jerri1977,Benn-samp,zayedbook,Marks2012}. Our aim here is to use Shannon sampling methods to map the discrete structure given by the primes into continuous functions which can then be Fourier analyzed.  

Concretely, our study has three parts. In the first part, we consider a histogram of the spacings between any pair of primes within some finite interval. In the second part, we use the primes to construct a continuous function by using Shannon sampling theory which is then Fourier analyzed. In the third part, we apply a generalized Shannon sampling method. 

With each method, we find that the primorial's role as jumping champions is foreshadowed in the sense that Fourier spikes for the primorials arise much earlier on the number axis than where the primorials in question are expected to reign as jumping champion. In addition, we also find prominent Fourier spikes at frequencies that would correspond to certain non-integer spacings. These spacing are simple ratios that are as yet unexplained.   

\section{Histogram analysis of Prime Gaps} \label{sec:analysishisto}
\begin{figure}[H]
    \includegraphics[width = \textwidth]{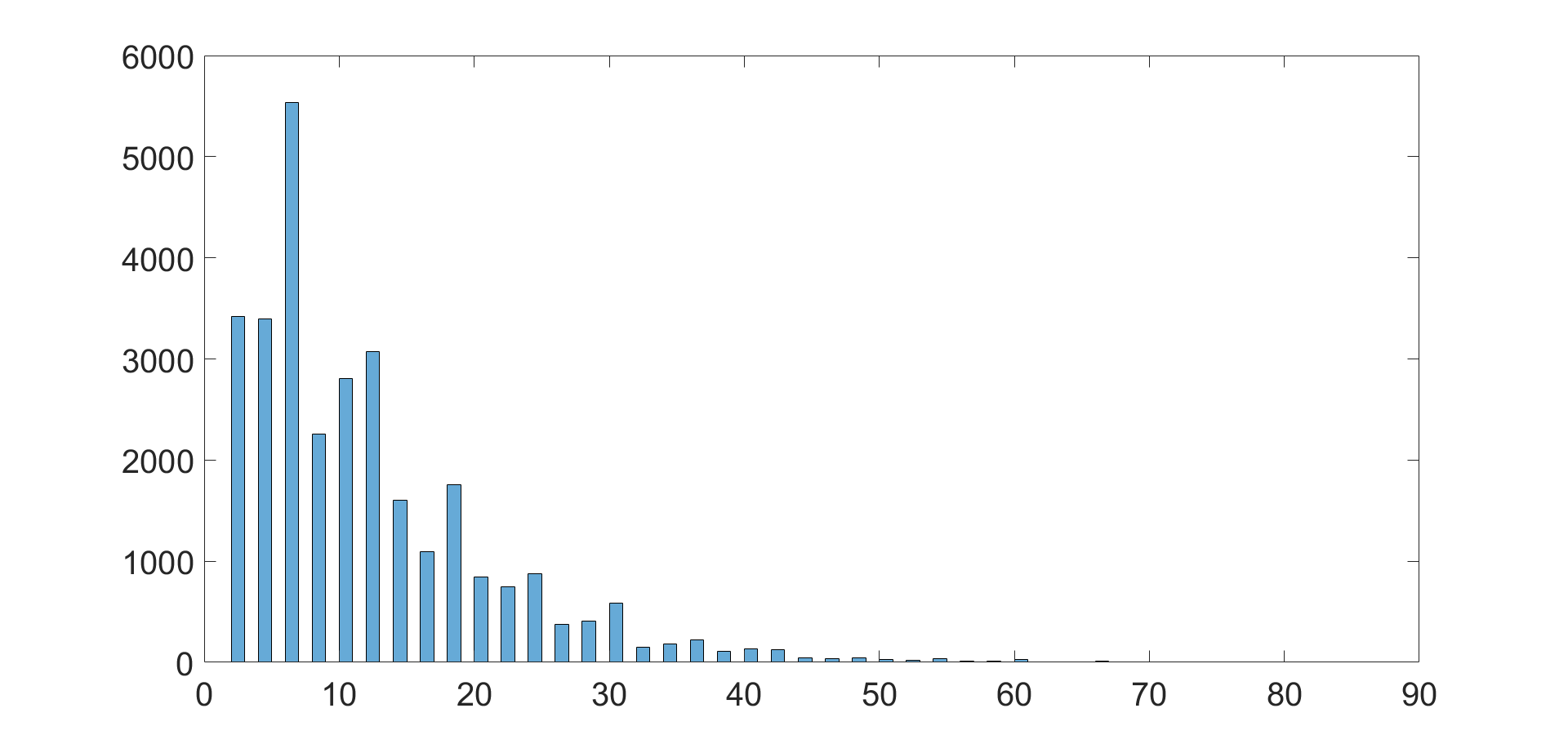}
    \caption{Distribution of the differences between consecutive primes, up to the 50000th prime} 
    \label{fig:hist1}
\end{figure}

We begin our analysis of the prime gaps by plotting a histogram of the differences between consecutive primes up to a maximum prime (Figure \ref{fig:hist1}).

\begin{figure}[H]
    \includegraphics[width = \textwidth]{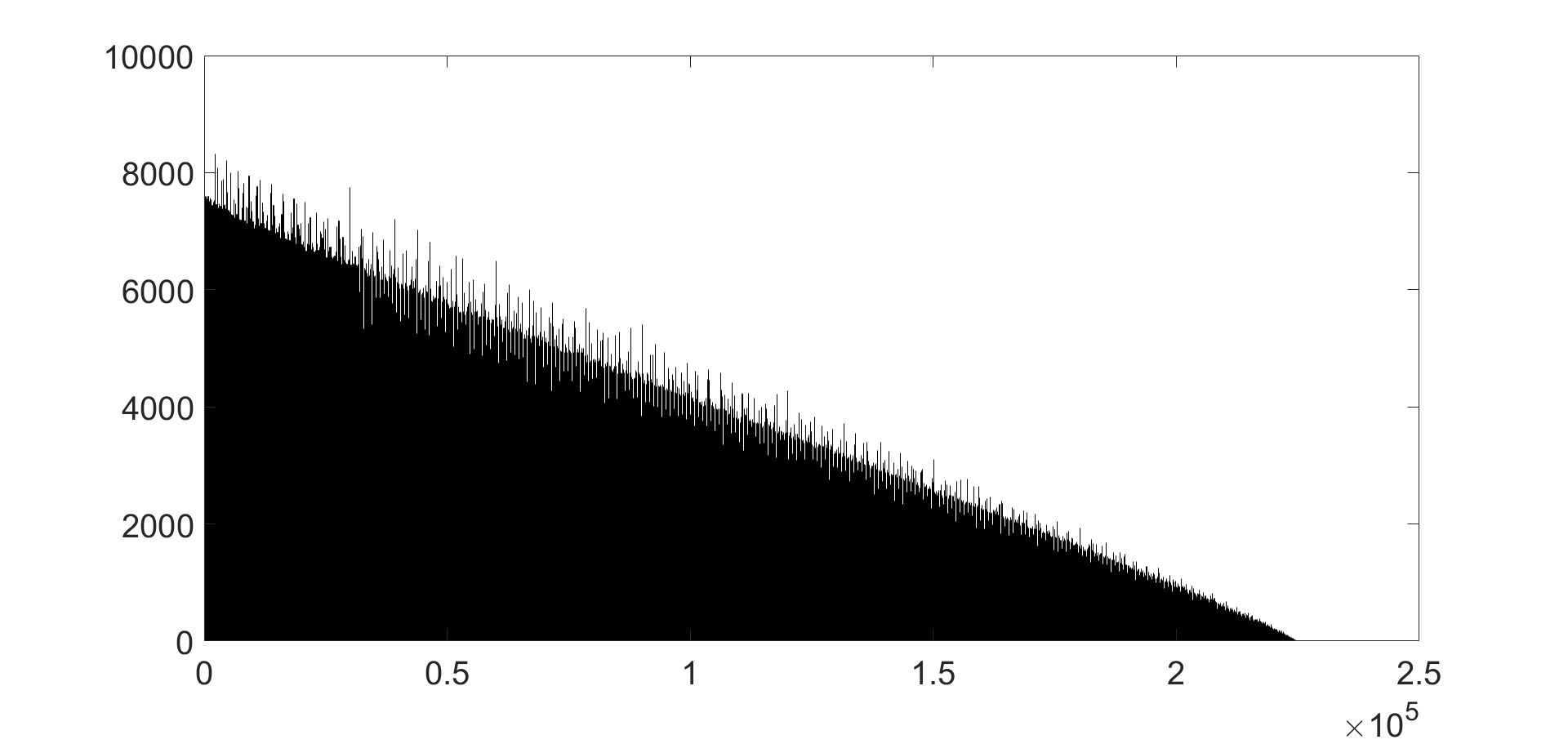}
    \caption{Distribution of the differences between every combination of primes, beginning from the second prime ($p_2 = 3$) to the 50000th prime. We find that some of the prominent spikes are at distance values that ``foreshadow" the occurrence of jumping champions.} 
    \label{fig:hist2}
\end{figure}

Previously, Ares and Castro \cite{aresandcastro} noticed that periodic oscillations occur within the histogram, leading to spikes at differences that are a multiple of 6. 
Motivated by this observation, let us now examine if similar structural properties exist among the differences between primes that need not be consecutive, see (Figure \ref{fig:hist2}).
Figure \ref{fig:hist2} exhibits a significant number of spikes and a closer inspection shows that among them are all the primorials up until the primorial 2310. This very early occurrence of such large primorials suggest closer examination with more powerful tools. To this end, we will now employ Shannon sampling methods from information-theory. We will use both regular and adaptive Shannon methods to map the discrete distribution of the primes into continuous functions whose periodicities, which are related to spacings of primes, can then be analyzed using the Fourier transform. 

\section{Frequency Analysis using Classical Shannon Sampling Theory} \label{sec:analysisclassical}
In information theory, the Shannon sampling theorem plays a central role as it establishes an equivalence of discrete and continuous representations of information. Concretely, it allows one to perfectly reconstruct a bandlimited (and therefore continuous) function on the real line from knowledge of its amplitudes only at a discrete set of points on the real line. We recall that a function over the reals is called bandlimited if the support of its Fourier transform is bounded.
We will use the Shannon sampling theorem to construct a continuous and bandlimited function by specifying its samples on the integers to be either 1 or 0 depending on whether the integer is prime or not. We then Fourier analyze the constructed function. 

\subsection{Background} 
The Shannon sampling theorem \cite{shannon,Jerri1977,Benn-samp,zayedbook,Marks2012} states that if a function, or `signal', $\phi(t)$, possesses no frequencies above some finite value $\Omega_{max}$, then it suffices to record the samples $\{\phi(t_n)\}$ on an equidistantly-spaced lattice $\{t_n\}$ with spacing $t_{n+1}-t_n=(2\Omega_{max})^{-1}$. If the samples are taken at this rate, the so-called Nyquist Rate, the function $\phi(t)$ can be reconstructed for all $t$ through
\begin{equation} \label{eq:recon}
\phi(t) = \sum_{n = -\infty}^{\infty}G(t, t_n)\phi(t_n)
\end{equation}
where $G(t, t_n)$, the reconstruction kernel, is defined as:
\begin{equation}
G(t, t_n) = \sinc\left(2(t - t_n)\Omega_{max}\right)
\end{equation}
We remark that $\Omega_{max}$-bandlimited functions can also be reconstructed from non-equidistantly spaced samples, if their average density (in the sense of Beurling) matches or exceeds the Nyquist density. The reconstruction from a non-equidistantly-spaced sampling lattice is necessarily less stable, however, in the sense that small measurement errors in the amplitudes translate into increased errors in the reconstructed function.


\subsection{Signal Construction and Analysis}
\begin{figure}[H]
    \includegraphics[width = \textwidth]{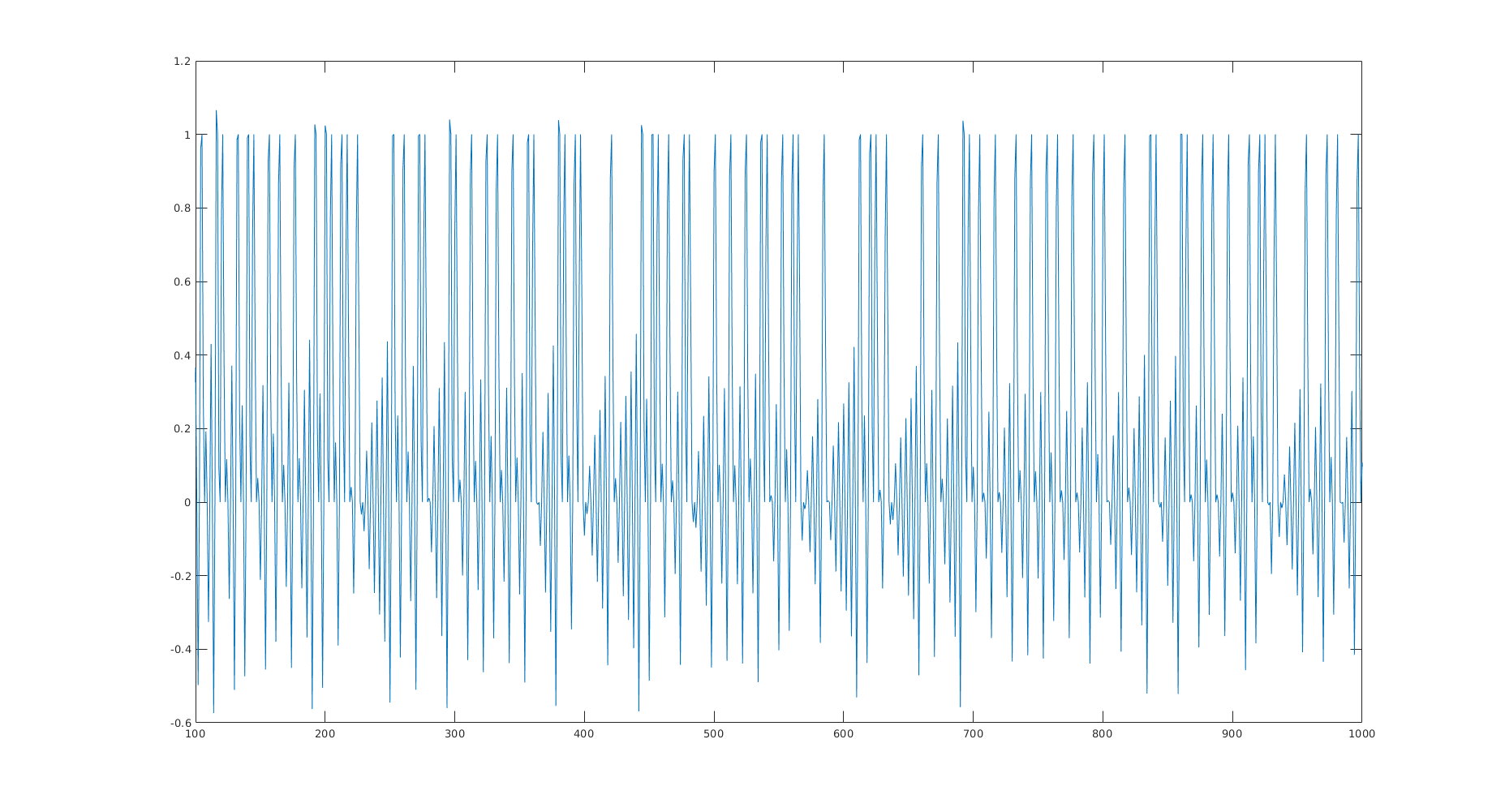}
    \caption{Zoom-In of prime number signal generated using the Shannon sampling theorem}
    \label{fig:reg_shannon}
\end{figure}

\begin{figure}[H]
    \includegraphics[width = \textwidth]{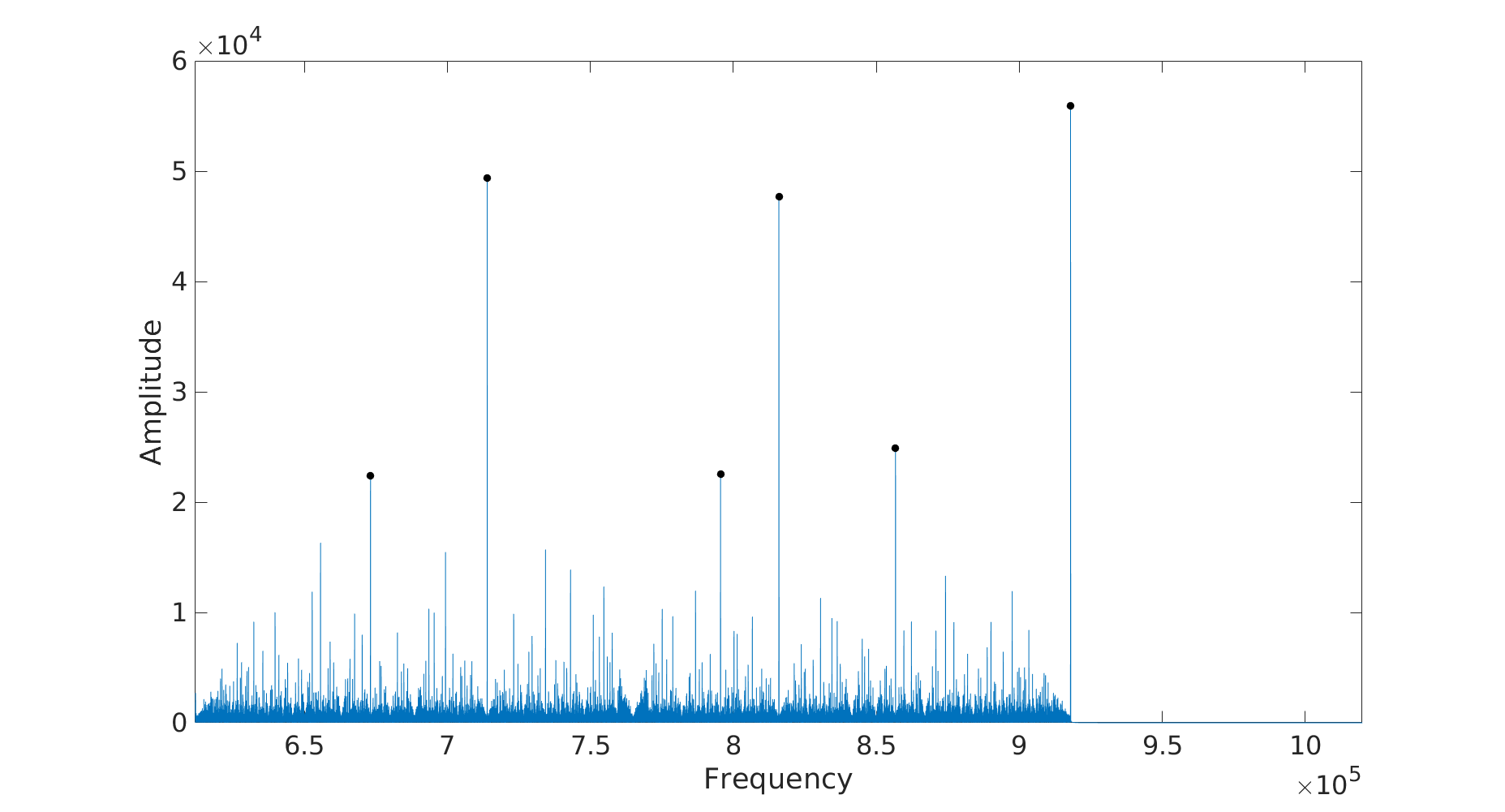}
    \caption{Fourier transform of Prime Signal Generated Using Classical Reconstruction Methods.Only the positive frequencies are shown, and the leftmost frequencies correspond to longer wavelengths. The most prominent spikes are dotted.}
    \label{fig:fourier_reg_shannon}
\end{figure}

Our aim is to construct a continuous function, or `signal', $\Phi(t)$, based on the primes. To this end, we define our sampling lattice $\{t_n\}$ to be the set of integers and we define the function's amplitudes on the integers to be:  
\begin{gather*}
\phi(t_n) =
\begin{cases}
0, \text{if $t_n$ non-prime} \\
1, \text{if $t_n$ prime}
\end{cases}
\end{gather*}
The resulting signal obtained by using the first 50000 primes is shown in Figure \ref{fig:reg_shannon} and its Fourier spectrum is shown in Figure \ref{fig:fourier_reg_shannon}. 

\begin{table}[H]
    \centering
       \begin{tabular}{c||c|c|c} 
        No. & Frequency (Hz) & Amplitude & Wavelength \\ 
        \hline
        1 & 673148 & 4941.04 & 20.000 \\ 
        2 & 713945 & 2241.09 & 12.000 \\ 
        3 & 734343 & 2252.28 & 6.6667 \\ 
        4 & 795539 & 4766.72 & 6.0000 \\ 
        5 & 856734 & 2490.08 & 5.0000\\ 
        6 & 917929 & 5591.44 & 4.0000  \\ 
    \end{tabular}
    \label{tab:regulargaps}
    \caption{The frequency, amplitude and wavelength of the most prominent spikes, from left to right in Figure \ref{fig:fourier_reg_shannon}.}
\end{table}
In Fig.\ref{fig:fourier_reg_shannon} and Table 1, the occurrence of spikes at certain integer wavelengths indicates the prevalence of corresponding prime gaps. In particular, 6 appears among those with highest amplitude, indicating that it is one of the most commonly occurring prime gaps, consistent with the observations made in Section \ref{sec:analysishisto}. We notice also that there are spikes at non-integer wavelengths which appear to be simple ratios. These non-integer `effective gaps' indicate the existence of structures that cannot be seen in a histogram of integer prime distances.

\section{Frequency analysis using generalized, adaptive Shannon sampling theory}
Our aim now is to try to further amplify 
the phenomenon of the foreshadowing of jumping champions and the occurrence of non-integer prominent prime spacings
by applying a generalized Shannon method that allows one to adapt the choice of sample points to arbitrary irregular spacings. This allows us to choose our sample points to be the sequence of the primes. We will again set the amplitude at the sample points to be 1 to obtain a signal that can then be Fourier analyzed. Notice that if we also require the amplitude to be 0 at the non-prime integers, then we recover our previous signal exactly. So, our new method is distinguished from the method of the previous section by in this sense only focusing on the prime numbers.
\label{sec:analysisgensamp}
\subsection{Background}
The generalized sampling theory \cite{Kempf2000,Kempf-samp,Hao2,yk1, yk2, yk3,KM,yufangthesis} generalizes the regular Shannon sampling theorem of functions, or signals, that possess a constant bandwidth and constant Nyquist rate to classes of functions that possess a time-varying bandwidth, or, correspondingly, a time-varying Nyquist rate. 
This allows one to consider classes of signals of time-varying bandwidth that can be most stably reconstructed from their amplitudes on a sampling lattice, $\{t_n\}$, whose spacing correspondingly varies in time. 
In addition to specifying the sampling lattice $\{t_n\}$, which we will choose to be the primes, the generalized sampling theory also requires the specification of 
a set of values $\{t_n'\}$ which in effect describe the extent to which the bandwidth may change from sample to sample. In the absence of additional information that we could use here, we will set these values to be the standard values of $t'_n=t_{n+1}-t_n$.

The reconstruction formula Eq. (\ref{eq:recon}) can now be applied with the generalized reconstruction kernel\cite{Kempf2000,Kempf-samp,Hao2,yk1, yk2, yk3,KM,yufangthesis}:
\begin{equation}
G(t, t_n) = (-1)^{z(t, t_n)}
\frac{\sqrt{t'_n}}{t - t_n}
\left(
\sum_m \frac{t'_m}{(t-t_m)^2}
\right)^2
\end{equation}
The function $z(t,t_n)$, in the exponent is the number of sampling points between $t$ and $t_n$, so that $(-1)^z(t, t_n)$ makes $G(t, t_n)$ differentiable.

More generally, given $s,t \in \mathbb{R}$, 
$$ G (s,t) := f(t) \cdot \left( \sum \frac{t_n '}{(t-t_n) (s-t_n)} \right) \cdot f(s), $$ where 
$$ f(t) := (-1) ^{z(t,t_n)} g(t) ^{-1/2}; \quad \quad g(t) := \sum \frac{t_n '}{(t-t_n) ^2}, $$ is a smooth (infinitely differentiable) positive kernel function on $\R \times \R$ in the sense of reproducing kernel Hilbert space (RKHS) theory, and our space of functions obeying a `generalized' bandlimit is the unique RKHS $\mc{H} (G)$ corresponding to $G$, see \cite[Section 2]{KM}. The theory of these spaces is closely connected to the theory of Hardy spaces of analytic functions in the complex upper-half plane \cite{KM,AMR,Clark1972,Krein1944,Krein1944one,Silva2007,Martin-dB,dB,GG,GR-model,Hoff}: Given any generalized RKHS, $\mc{H} (G)$, of time-varying or locally bandlimited functions, one can find a fixed function, $M(t)$, so that multiplication by $M(t)$ is a unitary transformation of $\mc{H} (G)$ onto a co-invariant (for the shift) subspace of the Hardy space of the upper half-plane which is the orthogonal complement of the range of a meromorphic inner function \cite{KM}.  This inner function can be expressed explicitly in terms of the sequences of sample points $\{ t_n \}$, and their `speeds' $\{ t_n '\}$ \cite{KM}.

The classical Shannon sampling kernel can be recovered as a special case of the generalized kernel $G(s,t)$ with the choice of sampling sequences:
$$ t_n := \frac{n\pi}{A}, \quad \mbox{and} \quad t_n ' := \tanh (A) \frac{\pi}{A}; \quad n \in \Z. $$ Indeed, one can apply trigonometric series identities to show
$$ \frac{1}{\pi} \sum _{k \in \Z} \left( \frac{1}{t-k} - \frac{1}{s-k} \right) = \cot(\pi t) - \cot(\pi s), $$ 
$$ g(t) ^{-2} = A \pi \tanh (A) \csc ^2 (At), $$ and finally that 
$$ G(s,t) = \frac{\sin \left( (t-s) A \right) }{(t-s) A}, $$ see \cite[Section 8.2]{Kempf-samp} or \cite[Example 2.28]{KM}.  

\subsection{Signal Construction and Analysis}
We construct our signal $\phi(t)$ as follows: our sampling lattice $\{t_n\}$ is a finite set of consecutive prime numbers, and $\phi(t_i) = 1$ for all $i$. Further, we set $\{t_n'\} = \{ \, (t_{i+1} - t_{i-1})/2 \, | \, i = 2, \ldots, n - 1 \, \}$ and $t_1' = (t_2 - t_1)/2$ and $t_n' = (t_n - t_{n-1})/2$.
As with our previous signal, we choose our sampling lattice to be the first 50000 primes, and the results of our reconstructed signal can be seen below. 
\begin{figure}[H]
    \includegraphics[width =\textwidth]{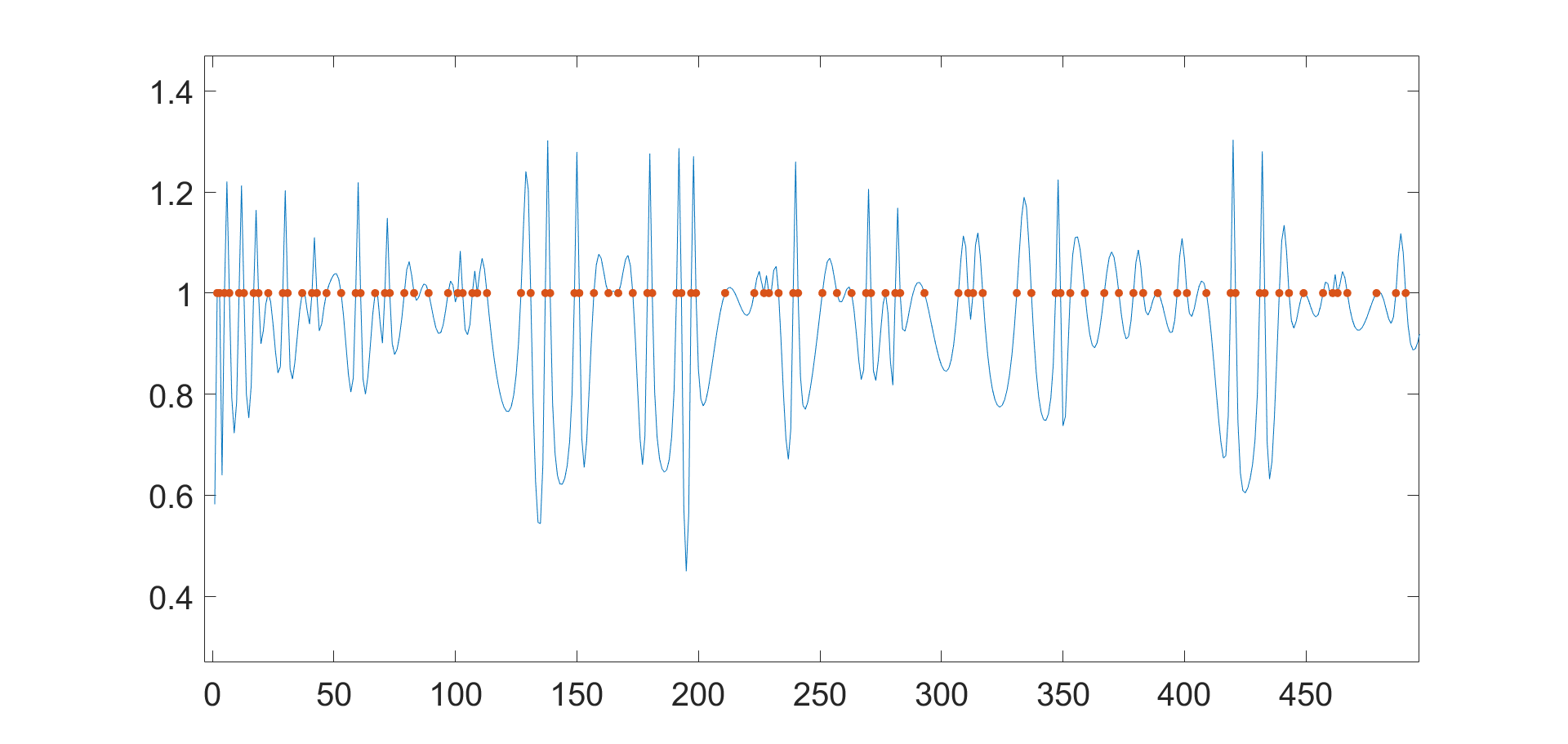}
    \caption{Zoom-in of the prime signal containing 50000 primes. The dotted points highlight where the prime numbers lie in the signal.}
    \label{fig:sig_section}
\end{figure}

\begin{figure}[H]
    \includegraphics[width = \textwidth]{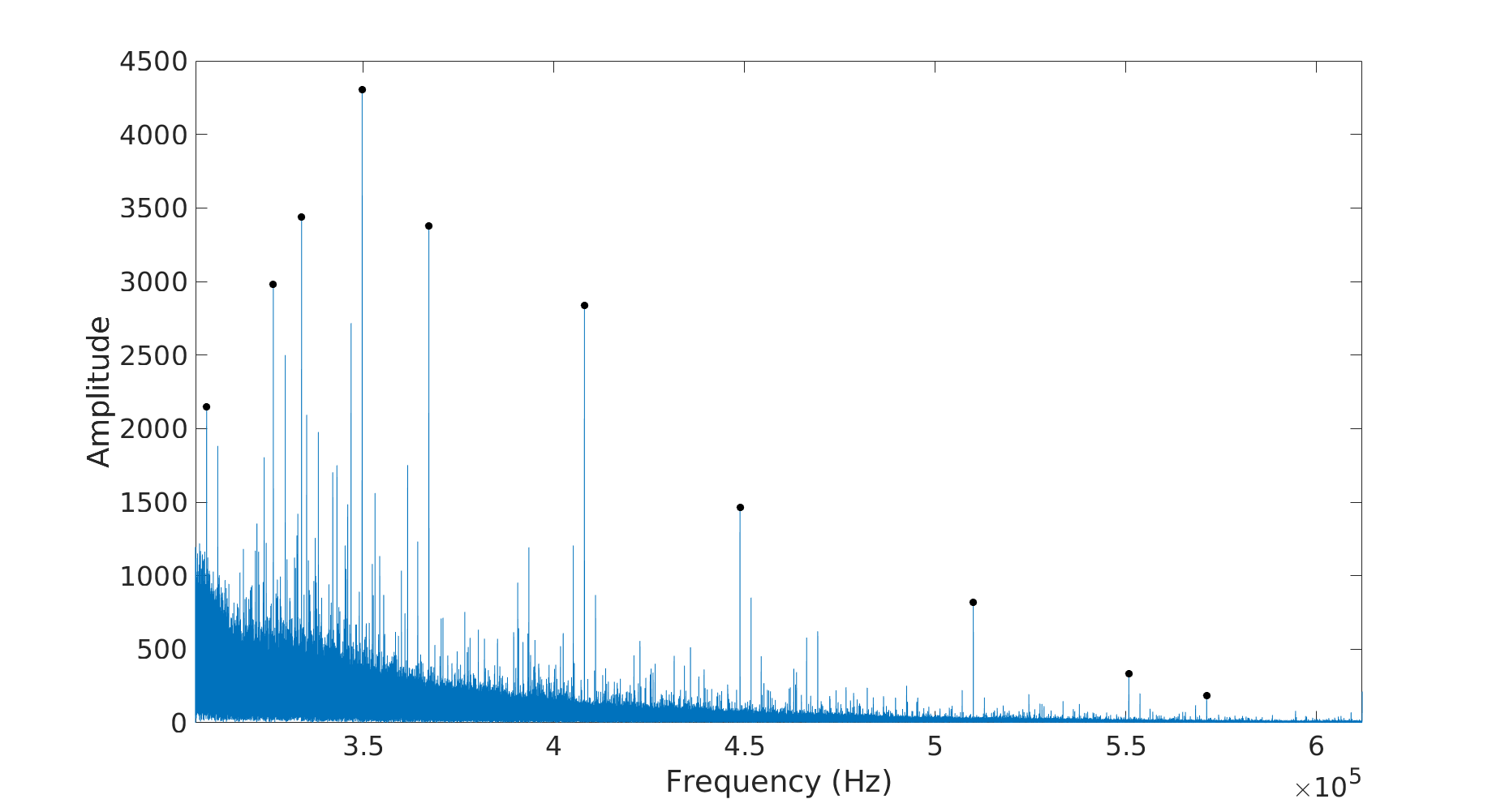}
    \caption{Modulus of the Fourier transform of the signal in Figure \ref{fig:sig_section}.}
    \label{fig:fft}
\end{figure}

\begin{table}[H] 
    \centering
    \begin{tabular}{c||c|c|c} 
        No. & Frequency (Hz) & Amplitude & Wavelength \\
        \hline
        1 & 2915 & 2146 & 210.02 \\ 
        2 & 20399 & 2981 & 30.001 \\ 
        3 & 27817 & 3438 & 22.000 \\ 
        4 & 43712 & 4306 & 14.000 \\ 
        5 & 61196 & 3377 & 10.000 \\ 
        6 & 101993 & 2837 & 6.0000 \\
        7 & 142790 & 1466 & 4.2857 \\ 
        8 & 203985 & 816.9 & 3.0000 \\ 
        9 & 244782 & 336.2 & 2.5000 \\ 
        10 & 265181 & 186.69 & 2.3077 \\
    \end{tabular}
    \caption{Calculated wavelengths from the Fourier transform in Figure \ref{fig:fft}. Notice that the non-integer spacing 4.2857 can be written as 30/7.} \label{tab:gaps1}
\end{table}
The result of using generalized sampling theory, compared to the previous two methods, is that we have uncovered the same and even more structure in the prime gaps. Previous results, such as the common gaps at multiples of 6, or the non-integer values are reproduced using this method. In addition, we find a foreshadowing of the next two conjectured jumping champions, 30 and 210. Given Conjecture \ref{conj:jc}), this suggests that the generalized sampling method is indeed able to accurately highlight frequent prime gaps. 

Let us now recall that in 2012, Goldon and Ledoan \cite{jumpingchampions2} provided several intervals in which a primorial is likely to become a jumping champion. The intervals are shown in Table \ref{tab:jc_ints}, where $p_k^{\#}$ denotes the product of the first $k$ primes (the $k$-th primorial). 
\begin{table}[H] 
    \centering
    \begin{tabular}{c||c|c|c} 
        No. & k & $p_k^{\#}$ & Interval \\
        \hline
        1 & 2 & 6 & $[4.67*10^4, 2.32*10^8]$\\ 
        2 & 3 & 30 & $[2.06 * 10^{44}, 5.24 * 10^{150}]$ \\ 
        3 & 4 & 210 & $[4.64 * 10^{487}, 4.01 * 10^{2607}]$\\ 
        4 & 5 & 2310 & $[8.78 * 10^{7769}, 1.72 * 10^{60178}]$ \\ 
        5 & 6 & 30030 & $[9.70 * 10^{134460}, 1.72 * 10^{1386286}]$ \\ 
    \end{tabular}
    \caption{Intervals in which $p_k^\#$ is likely to be a Jumping Champion} \label{tab:jc_ints}
\end{table}
These numbers indicate that in order to  confirm 30 as a jumping champion it is necessary to study the prime numbers up to about $10^{44}$. Yet, by creating a signal from only 50000 primes, a spike corresponding to the jumping champion 30 can already be observed in the Fourier transform. Likewise, seeing 210 foreshadowed as a spike is surprising. 

A possible explanation for the early occurrence of the primorials is that the primorials may not only be the jumping champions of subsequent primes but may also be closely related to the most frequently occurring gaps between any pairs of primes.

\section{Stability Analysis}

\subsection{Comparison of the distribution of primes to a Poisson-distributed sequence of integers}

\begin{figure}[H]
    \includegraphics[width = \textwidth]{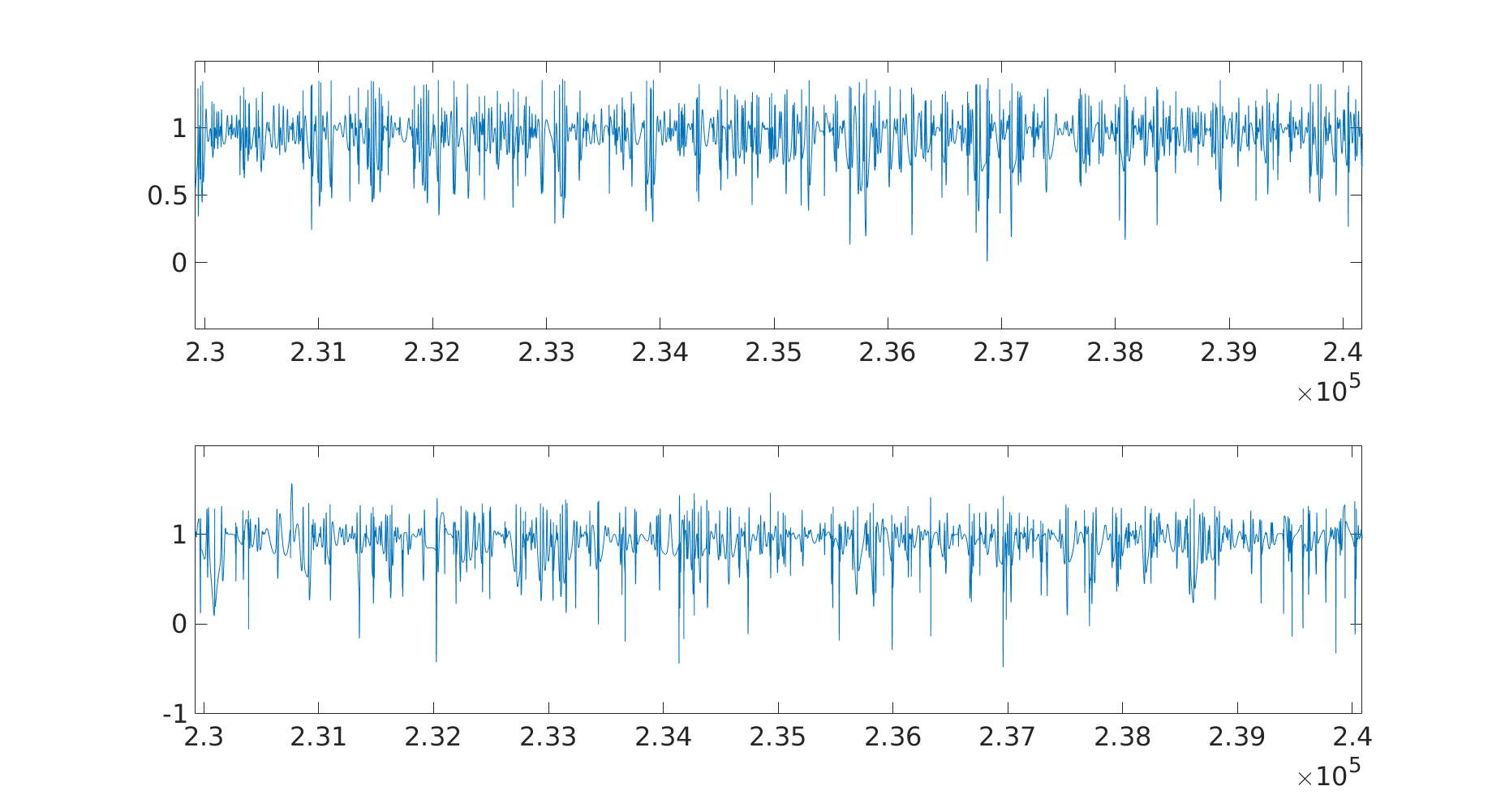}
    \caption{Top: Section of signal constructed from the first 50000 primes, Bottom: Section of signal constructed from 50000 Poisson-distributed numbers}
    \label{fig:poisson-signal}
\end{figure}

\begin{figure}[H]
    \includegraphics[width = \textwidth]{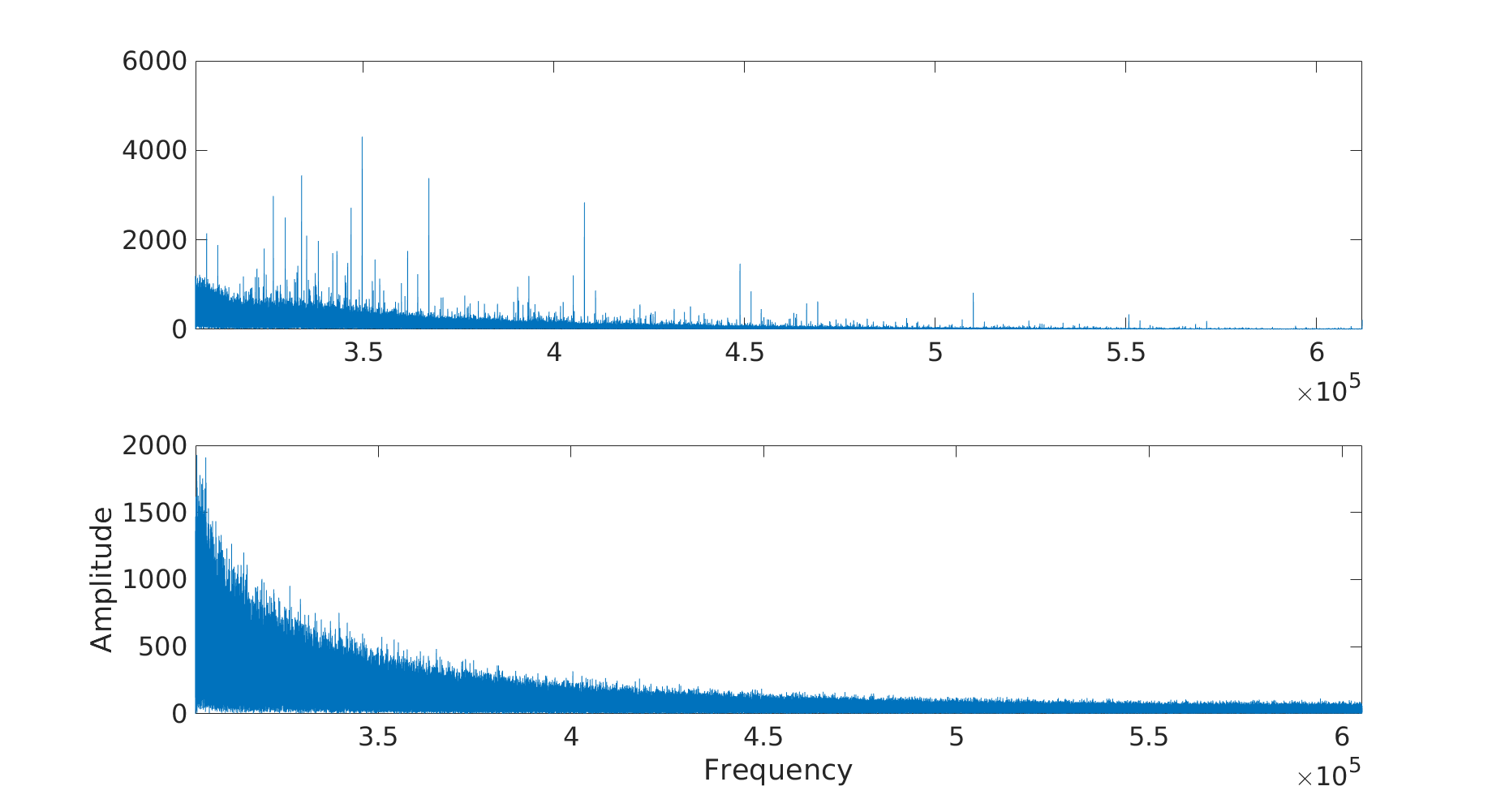}
    \caption{Top: Fourier transform of prime signal, Bottom: Fourier transform of poisson-distributed signal}
    \label{fig:poisson-fourier}
\end{figure}

For comparison, we now apply our method also to Poisson-distributed sets of integers of the same density as the primes, see Figure \ref{fig:poisson-signal}. We observe in Figure \ref{fig:poisson-fourier} that there is an accumulation of low frequencies in the Poisson case which in the case of the prime numbers become re-distributed to become the prominent spikes. 

\subsection{Verification of the translation-invariance of the generalized sampling methods}

Here, we apply our methods to three sampling lattices from various sections of the prime numbers: the 1st to 10000th prime, the 10001st to 20000th prime, and the 20001st to 30000th prime. We then calculate the Fourier transforms of the signals constructed from these sampling lattices. We observe that the prominent spikes are stable in the sense that they are prominent in all three Fourier transforms.

\begin{figure}[H]
    \includegraphics[width = \textwidth]{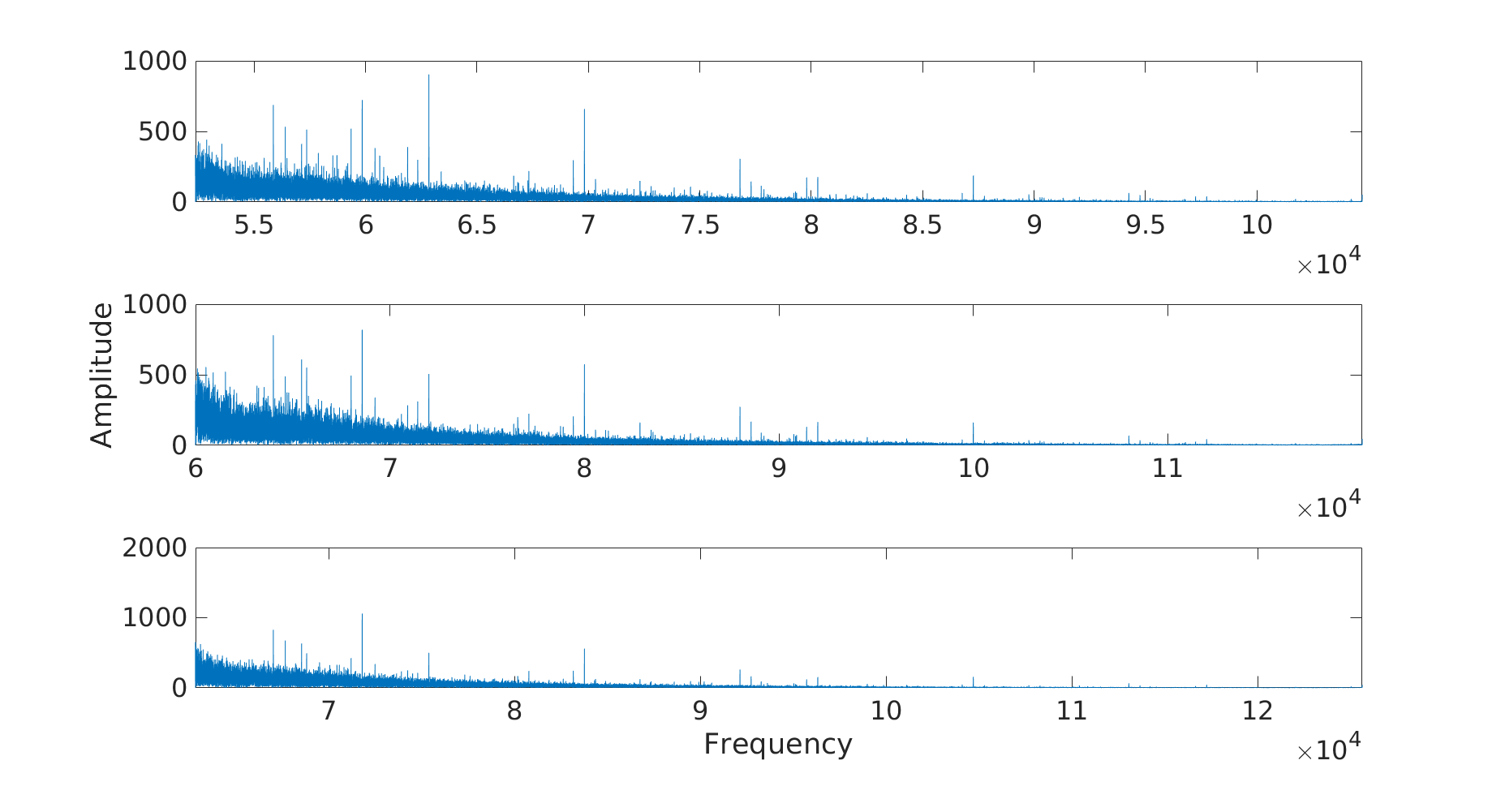}
    \caption{Top: Fourier transform of signal constructed from the 1st to 10000th prime \newline Center: Fourier transform of signal constructed from the 10001st to 20000th prime \newline Bottom: Fourier transform of signal constructed from the 20001st to 30000th prime}
    \label{fig:translation-invariance}
\end{figure}

\section{Conclusions and Outlook}

We analyzed the discrete distribution of the primes by Fourier analyzing continuous functions obtained from the primes. In order to map the sequence of primes into continuous functions, we used Shannon sampling methods from information theory. The local behavior of the information-theoretically obtained continuous function depends on the primes nonlocally, with the influence of primes a distance of $d$ away naturally decaying as $1/d$.

The Fourier transform of these continuous `prime signals' yielded intriguing peaks at the primorials, as well as at non-integer wavelengths. In particular, the application of adaptive Shannon methods yielded more and stronger peaks at the primorials. 

We conclude that, for as yet unknown reasons, the presence of the jumping champions of Conjecture \ref{conj:jc}, is foreshadowed far earlier on the number line, among very much smaller primes than expected.  This suggests that the primorials may also play a prominent role for the distances among any two not necessarily subsequent primes, giving rise to long-range correlations.

Also, intriguingly, our results show prominent wavelengths in the Fourier analysis of the prime signals that occur at values that are not integer and that therefore cannot directly correspond to prime gaps. These wavelengths, which may be called effective prime gaps, appear to be particularly simple ratios whose origin and structure should be very interesting to explore, as they may be related to the Chebychev bias in the distribution of primes, see, e.g. \cite{granville2006prime}, or more generally to the biases that were recently discovered in \cite{oliver2016unexpected}. 
\\ \\ \noindent
{\bf Application of generalized Shannon sampling method to other sequences} \\ \\
It should also be very interesting to apply the new method that uses adaptive Shannon sampling to other sequences. For example, we have applied the new method to the sequence of squares and twice the squares (SEQ1), and the sequence of integers that are the sums of two squares (SEQ2), as shown in the figures below.

\begin{figure}[H]
    \includegraphics[width = \textwidth]{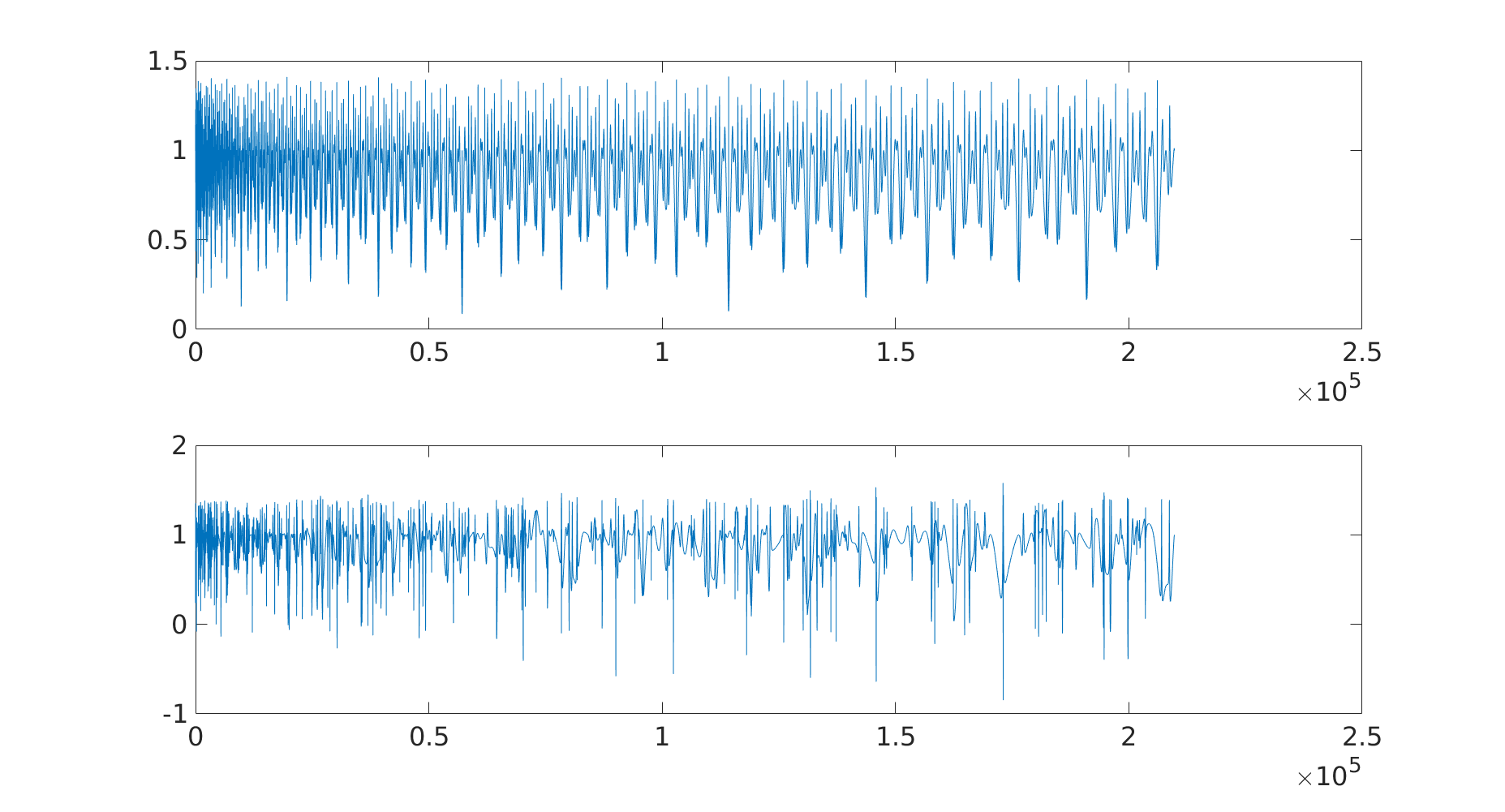}
    \caption{Top: Signal constructed of numbers from SEQ1 less than 209760. 
    Bottom: Signal constructed from Poisson-distributed numbers of the same density as SEQ1}
    \label{fig:sq+2sq-signal}
\end{figure}

\begin{figure}[H]
    \includegraphics[width = \textwidth]{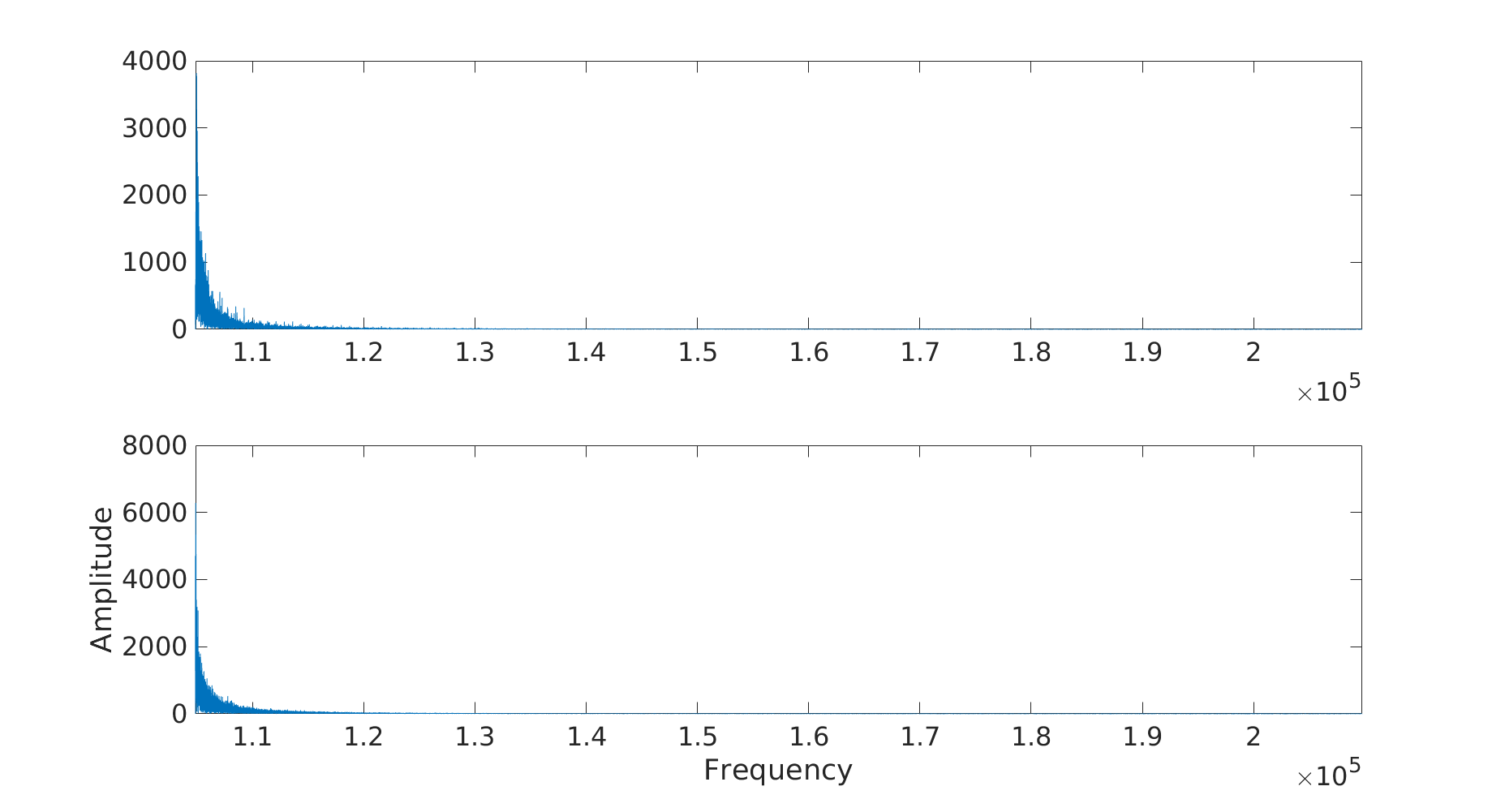}
    \caption{Top: Fourier transform of signal constructed from SEQ1.
    Bottom: Fourier transform of Poisson-distributed numbers of the same density as SEQ1}
    \label{fig:sq+2sq-fourier}
\end{figure}

\begin{figure}[H]
    \includegraphics[width = \textwidth]{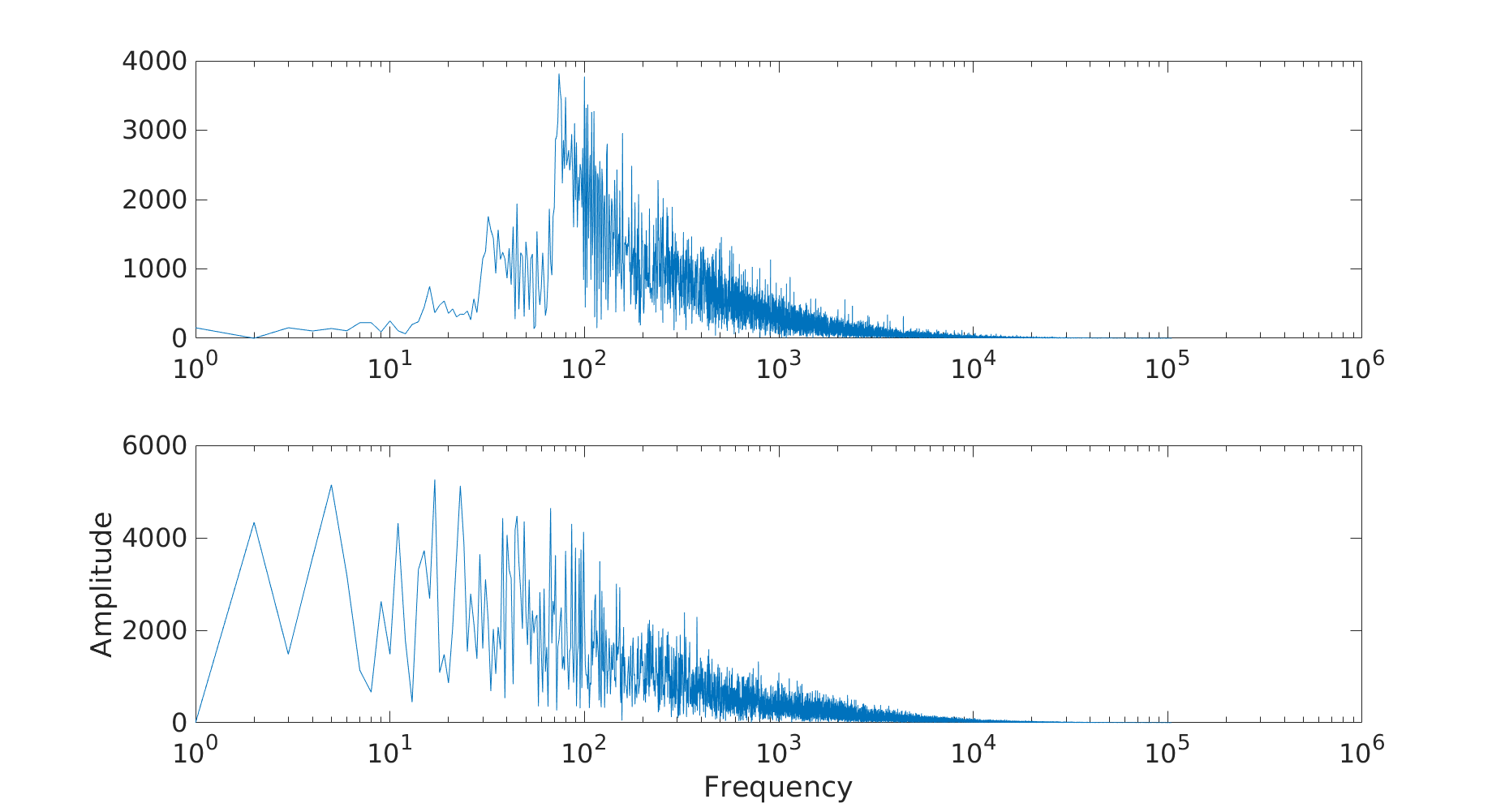}
    \caption{Top: Fourier transform of signal constructed from SEQ1 with a logarithmic scale for the x-axis
    Bottom: Fourier transform of Poisson-distributed numbers of the same density as SEQ1 with a logarithmic scale for the x-axis}
    \label{fig:sq+2sq-fourier-logscale}
\end{figure}

\begin{figure}[H]
    \includegraphics[width = \textwidth]{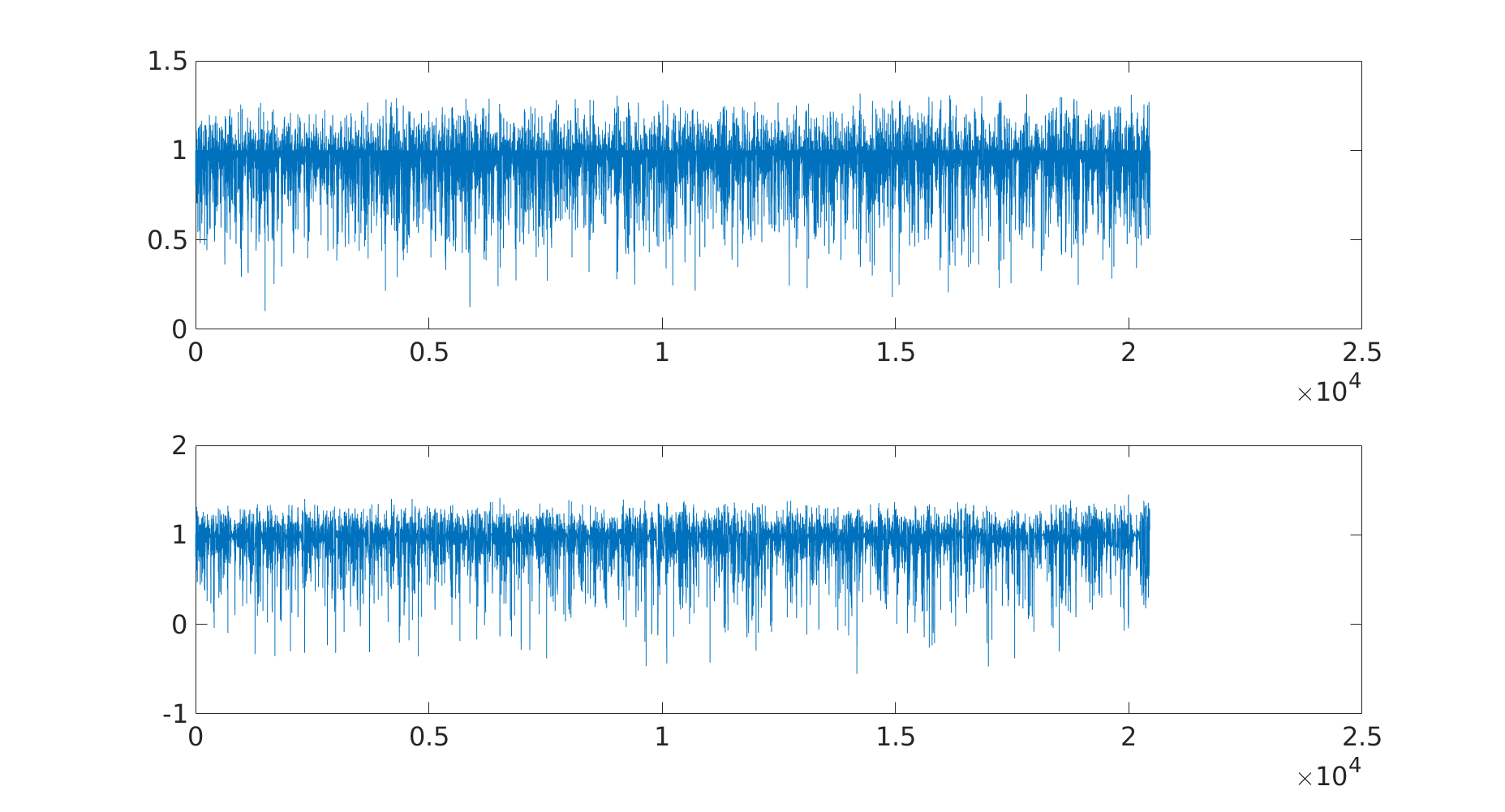}
    \caption{Top: Signal constructed from the first 5871 integers of SEQ2.
    Bottom: Signal constructed from Poisson-distributed numbers of the same density as SEQ2}
    \label{fig:sumsq-signal}
\end{figure}

\begin{figure}[H]
    \includegraphics[width = \textwidth]{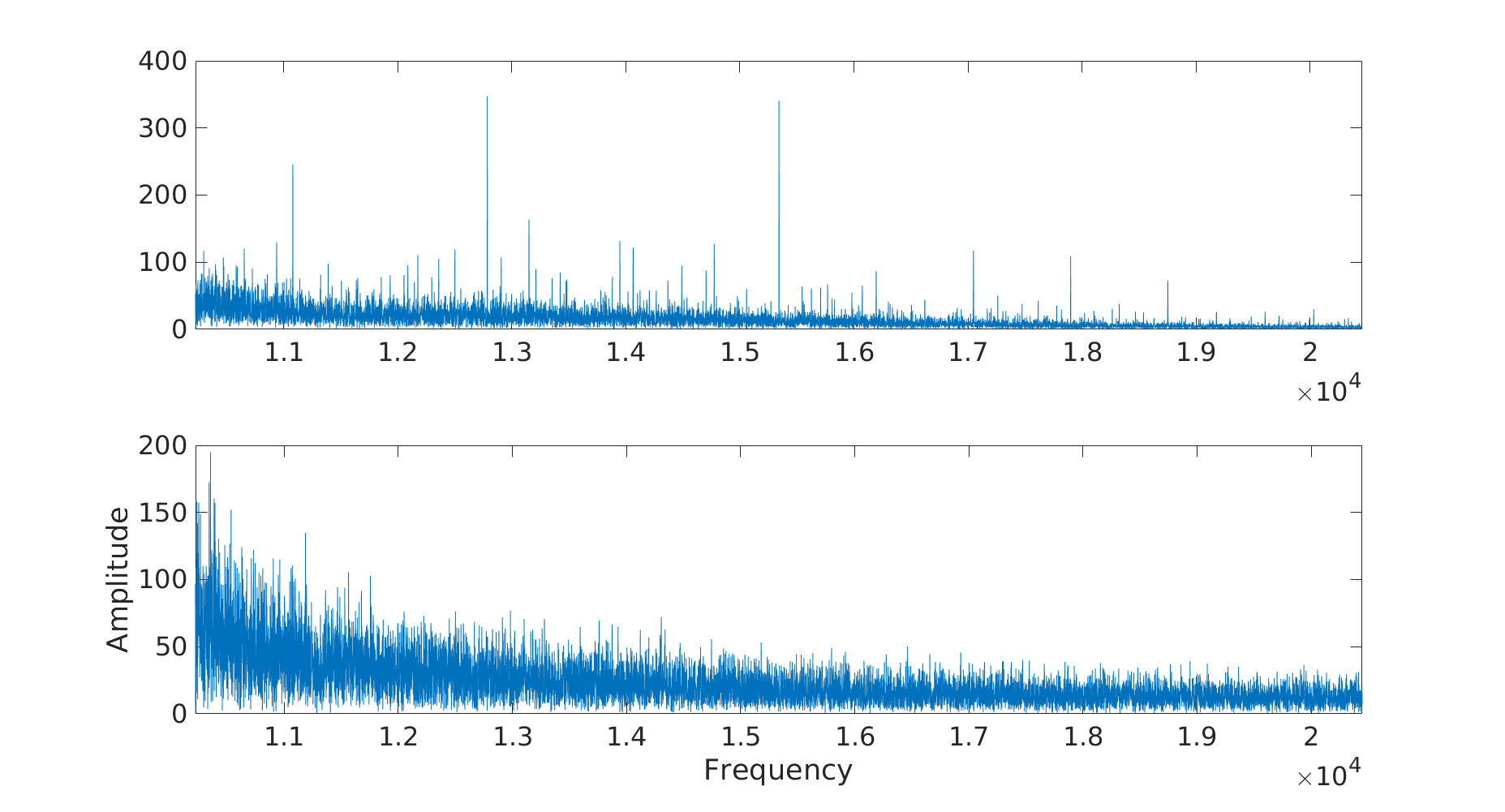}
    \caption{Top: Fourier transform of signal constructed from SEQ2.
    Bottom: Fourier transform of Poisson-distributed numbers of the same density as SEQ2}
    \label{fig:sumsq-fourier}
\end{figure}

$$$$
\bf Acknowledgements: \rm NP, AK and RTW are grateful for useful feedback from Kevin Hare, Tristan Freiberg and Stefan Steinerberger. AK acknowledges support from the Discovery program of the National Science and Engineering Research Council of Canada (NSERC). 

\section*{References}


\begin{thebibliography}{33}%
\makeatletter
\providecommand \@ifxundefined [1]{%
 \@ifx{#1\undefined}
}%
\providecommand \@ifnum [1]{%
 \ifnum #1\expandafter \@firstoftwo
 \else \expandafter \@secondoftwo
 \fi
}%
\providecommand \@ifx [1]{%
 \ifx #1\expandafter \@firstoftwo
 \else \expandafter \@secondoftwo
 \fi
}%
\providecommand \natexlab [1]{#1}%
\providecommand \enquote  [1]{``#1''}%
\providecommand \bibnamefont  [1]{#1}%
\providecommand \bibfnamefont [1]{#1}%
\providecommand \citenamefont [1]{#1}%
\providecommand \href@noop [0]{\@secondoftwo}%
\providecommand \href [0]{\begingroup \@sanitize@url \@href}%
\providecommand \@href[1]{\@@startlink{#1}\@@href}%
\providecommand \@@href[1]{\endgroup#1\@@endlink}%
\providecommand \@sanitize@url [0]{\catcode `\\12\catcode `\$12\catcode
  `\&12\catcode `\#12\catcode `\^12\catcode `\_12\catcode `\%12\relax}%
\providecommand \@@startlink[1]{}%
\providecommand \@@endlink[0]{}%
\providecommand \url  [0]{\begingroup\@sanitize@url \@url }%
\providecommand \@url [1]{\endgroup\@href {#1}{\urlprefix }}%
\providecommand \urlprefix  [0]{URL }%
\providecommand \Eprint [0]{\href }%
\providecommand \doibase [0]{http://dx.doi.org/}%
\providecommand \selectlanguage [0]{\@gobble}%
\providecommand \bibinfo  [0]{\@secondoftwo}%
\providecommand \bibfield  [0]{\@secondoftwo}%
\providecommand \translation [1]{[#1]}%
\providecommand \BibitemOpen [0]{}%
\providecommand \bibitemStop [0]{}%
\providecommand \bibitemNoStop [0]{.\EOS\space}%
\providecommand \EOS [0]{\spacefactor3000\relax}%
\providecommand \BibitemShut  [1]{\csname bibitem#1\endcsname}%
\let\auto@bib@innerbib\@empty
\bibitem [{\citenamefont {Hardy}\ and\ \citenamefont
  {Littlewood}(1923)}]{hardy1923}%
  \BibitemOpen
  \bibfield  {author} {\bibinfo {author} {\bibfnamefont {G.~H.}\ \bibnamefont
  {Hardy}}\ and\ \bibinfo {author} {\bibfnamefont {J.~E.}\ \bibnamefont
  {Littlewood}},\ }\bibfield  {title} {\enquote {\bibinfo {title} {Some
  problems of partitio numerorum; iii: On the expression of a number as a sum
  of primes.}}\ }\href@noop {} {\bibfield  {journal} {\bibinfo  {journal} {Acta
  Math.}\ }\textbf {\bibinfo {volume} {44}},\ \bibinfo {pages} {1--70}
  (\bibinfo {year} {1923})}\BibitemShut {NoStop}%
\bibitem [{\citenamefont {Odlyzko}, \citenamefont {Rubinstein},\ and\
  \citenamefont {Wolf}(1999)}]{jumpingchampions}%
  \BibitemOpen
  \bibfield  {author} {\bibinfo {author} {\bibfnamefont {A.}~\bibnamefont
  {Odlyzko}}, \bibinfo {author} {\bibfnamefont {M.}~\bibnamefont {Rubinstein}},
  \ and\ \bibinfo {author} {\bibfnamefont {M.}~\bibnamefont {Wolf}},\
  }\bibfield  {title} {\enquote {\bibinfo {title} {Jumping champions},}\
  }\href@noop {} {\bibfield  {journal} {\bibinfo  {journal} {Experimental
  Mathematics}\ }\textbf {\bibinfo {volume} {8}},\ \bibinfo {pages} {107--118}
  (\bibinfo {year} {1999})}\BibitemShut {NoStop}%
\bibitem [{\citenamefont {Goldston}\ and\ \citenamefont
  {Ledoan}(2015)}]{jumpingchampions2}%
  \BibitemOpen
  \bibfield  {author} {\bibinfo {author} {\bibfnamefont {D.}~\bibnamefont
  {Goldston}}\ and\ \bibinfo {author} {\bibfnamefont {A.}~\bibnamefont
  {Ledoan}},\ }\bibfield  {title} {\enquote {\bibinfo {title} {The jumping
  champion conjecture},}\ }\href@noop {} {\bibfield  {journal} {\bibinfo
  {journal} {Mathematika}\ }\textbf {\bibinfo {volume} {61}},\ \bibinfo {pages}
  {719--740} (\bibinfo {year} {2015})},\ \Eprint
  {http://arxiv.org/abs/1102.4879v2} {arXiv:1102.4879v2 [math.NT]} \BibitemShut
  {NoStop}%
\bibitem [{\citenamefont {Wolf}(1999)}]{Wolf1999}%
  \BibitemOpen
  \bibfield  {author} {\bibinfo {author} {\bibfnamefont {M.}~\bibnamefont
  {Wolf}},\ }\bibfield  {title} {\enquote {\bibinfo {title} {Applications of
  statistical mechanics in number theory},}\ }\href@noop {} {\bibfield
  {journal} {\bibinfo  {journal} {Physica A: Statistical Mechanics and its
  Applications}\ }\textbf {\bibinfo {volume} {274}},\ \bibinfo {pages}
  {149--157} (\bibinfo {year} {1999})}\BibitemShut {NoStop}%
\bibitem [{\citenamefont {Wolf}(2011)}]{Wolf2011}%
  \BibitemOpen
  \bibfield  {author} {\bibinfo {author} {\bibfnamefont {M.}~\bibnamefont
  {Wolf}},\ }\bibfield  {title} {\enquote {\bibinfo {title} {Some heuristics on
  the gaps between consecutive primes},}\ }\href@noop {} {\bibfield  {journal}
  {\bibinfo  {journal} {arXiv preprint arXiv:1102.0481}\ } (\bibinfo {year}
  {2011})}\BibitemShut {NoStop}%
\bibitem [{\citenamefont {Ares}\ and\ \citenamefont
  {Castro}(2006{\natexlab{a}})}]{Ares2006}%
  \BibitemOpen
  \bibfield  {author} {\bibinfo {author} {\bibfnamefont {S.}~\bibnamefont
  {Ares}}\ and\ \bibinfo {author} {\bibfnamefont {M.}~\bibnamefont {Castro}},\
  }\bibfield  {title} {\enquote {\bibinfo {title} {Hidden structure in the
  randomness of the prime number sequence?}}\ }\href@noop {} {\bibfield
  {journal} {\bibinfo  {journal} {Physica A: Statistical Mechanics and its
  Applications}\ }\textbf {\bibinfo {volume} {360}},\ \bibinfo {pages}
  {285--296} (\bibinfo {year} {2006}{\natexlab{a}})}\BibitemShut {NoStop}%
\bibitem [{\citenamefont {Szpiro}(2004)}]{Szpiro2004}%
  \BibitemOpen
  \bibfield  {author} {\bibinfo {author} {\bibfnamefont {G.}~\bibnamefont
  {Szpiro}},\ }\bibfield  {title} {\enquote {\bibinfo {title} {The gaps between
  the gaps: some patterns in the prime number sequence},}\ }\href@noop {}
  {\bibfield  {journal} {\bibinfo  {journal} {Physica A: Statistical Mechanics
  and its Applications}\ }\textbf {\bibinfo {volume} {341}},\ \bibinfo {pages}
  {607--617} (\bibinfo {year} {2004})}\BibitemShut {NoStop}%
\bibitem [{\citenamefont {Shannon}(1998)}]{shannon}%
  \BibitemOpen
  \bibfield  {author} {\bibinfo {author} {\bibfnamefont {C.~E.}\ \bibnamefont
  {Shannon}},\ }\bibfield  {title} {\enquote {\bibinfo {title} {Communication
  in the presence of noise},}\ }\href@noop {} {\bibfield  {journal} {\bibinfo
  {journal} {Proc. IEEE}\ }\textbf {\bibinfo {volume} {86}},\ \bibinfo {pages}
  {447--457} (\bibinfo {year} {1998})}\BibitemShut {NoStop}%
\bibitem [{\citenamefont {Jerri}(1977)}]{Jerri1977}%
  \BibitemOpen
  \bibfield  {author} {\bibinfo {author} {\bibfnamefont {A.}~\bibnamefont
  {Jerri}},\ }\bibfield  {title} {\enquote {\bibinfo {title} {The {S}hannon
  sampling theorem - its various extensions and applications: A tutorial
  review},}\ }\href@noop {} {\bibfield  {journal} {\bibinfo  {journal} {Proc.
  IEEE}\ }\textbf {\bibinfo {volume} {65}},\ \bibinfo {pages} {1565--1596}
  (\bibinfo {year} {1977})}\BibitemShut {NoStop}%
\bibitem [{\citenamefont {Benedetto}\ and\ \citenamefont
  {Ferreira}(2012)}]{Benn-samp}%
  \BibitemOpen
  \bibfield  {author} {\bibinfo {author} {\bibfnamefont {J.}~\bibnamefont
  {Benedetto}}\ and\ \bibinfo {author} {\bibfnamefont {P.}~\bibnamefont
  {Ferreira}},\ }\href@noop {} {\emph {\bibinfo {title} {Modern sampling
  theory: mathematics and applications}}}\ (\bibinfo  {publisher} {Springer},\
  \bibinfo {year} {2012})\BibitemShut {NoStop}%
\bibitem [{\citenamefont {Zayed}(1993)}]{zayedbook}%
  \BibitemOpen
  \bibfield  {author} {\bibinfo {author} {\bibfnamefont {A.~I.}\ \bibnamefont
  {Zayed}},\ }\href@noop {} {\emph {\bibinfo {title} {Advances in Shannon's
  Sampling Theory}}}\ (\bibinfo  {publisher} {CRC Press, Inc.},\ \bibinfo
  {year} {1993})\BibitemShut {NoStop}%
\bibitem [{\citenamefont {Marks}(2012)}]{Marks2012}%
  \BibitemOpen
  \bibfield  {author} {\bibinfo {author} {\bibfnamefont {R.}~\bibnamefont
  {Marks}},\ }\href@noop {} {\emph {\bibinfo {title} {Introduction to {S}hannon
  sampling and interpolation theory}}}\ (\bibinfo  {publisher} {Springer},\
  \bibinfo {year} {2012})\BibitemShut {NoStop}%
\bibitem [{\citenamefont {Ares}\ and\ \citenamefont
  {Castro}(2006{\natexlab{b}})}]{aresandcastro}%
  \BibitemOpen
  \bibfield  {author} {\bibinfo {author} {\bibfnamefont {S.}~\bibnamefont
  {Ares}}\ and\ \bibinfo {author} {\bibfnamefont {M.}~\bibnamefont {Castro}},\
  }\bibfield  {title} {\enquote {\bibinfo {title} {Hidden structure in the
  randomness of the prime number sequence?}}\ }\href@noop {} {\bibfield
  {journal} {\bibinfo  {journal} {Physica A}\ }\textbf {\bibinfo {volume}
  {360}},\ \bibinfo {pages} {285--296} (\bibinfo {year}
  {2006}{\natexlab{b}})},\ \Eprint {http://arxiv.org/abs/0310148v2}
  {arXiv:0310148v2 [cond-mat]} \BibitemShut {NoStop}%
\bibitem [{\citenamefont {Kempf}(2000)}]{Kempf2000}%
  \BibitemOpen
  \bibfield  {author} {\bibinfo {author} {\bibfnamefont {A.}~\bibnamefont
  {Kempf}},\ }\bibfield  {title} {\enquote {\bibinfo {title} {Fields over
  unsharp coordinates},}\ }\href@noop {} {\bibfield  {journal} {\bibinfo
  {journal} {Phys. Rev. Lett.}\ }\textbf {\bibinfo {volume} {85}},\ \bibinfo
  {pages} {2873} (\bibinfo {year} {2000})}\BibitemShut {NoStop}%
\bibitem [{\citenamefont {Kempf}(2004)}]{Kempf-samp}%
  \BibitemOpen
  \bibfield  {author} {\bibinfo {author} {\bibfnamefont {A.}~\bibnamefont
  {Kempf}},\ }\bibfield  {title} {\enquote {\bibinfo {title} {Fields with
  finite information density},}\ }\href@noop {} {\bibfield  {journal} {\bibinfo
   {journal} {Phys. Rev. D}\ }\textbf {\bibinfo {volume} {69}},\ \bibinfo
  {pages} {124014} (\bibinfo {year} {2004})}\BibitemShut {NoStop}%
\bibitem [{\citenamefont {Hao}\ and\ \citenamefont
  {Kempf}(2010{\natexlab{a}})}]{Hao2}%
  \BibitemOpen
  \bibfield  {author} {\bibinfo {author} {\bibfnamefont {Y.}~\bibnamefont
  {Hao}}\ and\ \bibinfo {author} {\bibfnamefont {A.}~\bibnamefont {Kempf}},\
  }\bibfield  {title} {\enquote {\bibinfo {title} {Generalized {S}hannon
  sampling method reduces the {G}ibbs overshoot in the approximation of a step
  function},}\ }\href@noop {} {\bibfield  {journal} {\bibinfo  {journal} {J.
  Concr. Appl. Math.}\ }\textbf {\bibinfo {volume} {8}},\ \bibinfo {pages}
  {540--554} (\bibinfo {year} {2010}{\natexlab{a}})}\BibitemShut {NoStop}%
\bibitem [{\citenamefont {Hao}\ and\ \citenamefont {Kempf}(2007)}]{yk1}%
  \BibitemOpen
  \bibfield  {author} {\bibinfo {author} {\bibfnamefont {Y.}~\bibnamefont
  {Hao}}\ and\ \bibinfo {author} {\bibfnamefont {A.}~\bibnamefont {Kempf}},\
  }\bibfield  {title} {\enquote {\bibinfo {title} {On a non-fourier
  generalization of shannon sampling theory},}\ }in\ \href@noop {} {\emph
  {\bibinfo {booktitle} {Information Theory, 2007. CWIT '07. 10th Canadian
  Workshop on}}}\ (\bibinfo {year} {2007})\ pp.\ \bibinfo {pages}
  {193--196}\BibitemShut {NoStop}%
\bibitem [{\citenamefont {Hao}\ and\ \citenamefont {Kempf}(2008)}]{yk2}%
  \BibitemOpen
  \bibfield  {author} {\bibinfo {author} {\bibfnamefont {Y.}~\bibnamefont
  {Hao}}\ and\ \bibinfo {author} {\bibfnamefont {A.}~\bibnamefont {Kempf}},\
  }\bibfield  {title} {\enquote {\bibinfo {title} {On the stability of a
  generalized shannon sampling theorem},}\ }in\ \href@noop {} {\emph {\bibinfo
  {booktitle} {Information Theory and Its Applications, 2008. ISITA 2008.
  International Symposium on}}}\ (\bibinfo {year} {2008})\ pp.\ \bibinfo
  {pages} {1--6}\BibitemShut {NoStop}%
\bibitem [{\citenamefont {Hao}\ and\ \citenamefont
  {Kempf}(2010{\natexlab{b}})}]{yk3}%
  \BibitemOpen
  \bibfield  {author} {\bibinfo {author} {\bibfnamefont {Y.}~\bibnamefont
  {Hao}}\ and\ \bibinfo {author} {\bibfnamefont {A.}~\bibnamefont {Kempf}},\
  }\bibfield  {title} {\enquote {\bibinfo {title} {Filtering, sampling, and
  reconstruction with time-varying bandwidths},}\ }\href@noop {} {\bibfield
  {journal} {\bibinfo  {journal} {IEEE Signal Proc. Lett.}\ }\textbf {\bibinfo
  {volume} {17}},\ \bibinfo {pages} {241--244} (\bibinfo {year}
  {2010}{\natexlab{b}})}\BibitemShut {NoStop}%
\bibitem [{\citenamefont {Martin}\ and\ \citenamefont {Kempf}(2018)}]{KM}%
  \BibitemOpen
  \bibfield  {author} {\bibinfo {author} {\bibfnamefont {R.}~\bibnamefont
  {Martin}}\ and\ \bibinfo {author} {\bibfnamefont {A.}~\bibnamefont {Kempf}},\
  }\bibfield  {title} {\enquote {\bibinfo {title} {Function spaces obeying a
  time-varying bandlimit},}\ }\href@noop {} {\bibfield  {journal} {\bibinfo
  {journal} {J. Math. Anal. Appl.}\ }\textbf {\bibinfo {volume} {458}},\
  \bibinfo {pages} {1597--1638} (\bibinfo {year} {2018})}\BibitemShut {NoStop}%
\bibitem [{\citenamefont {{Yufang Hao}}(2011)}]{yufangthesis}%
  \BibitemOpen
  \bibfield  {author} {\bibinfo {author} {\bibnamefont {{Yufang Hao}}},\ }\emph
  {\bibinfo {title} {{Generalizing Sampling Theory for Time-Varying Nyquist
  Rates using Self-Adjoint Extensions of Symmetric Operators with Deficiency
  Indices (1,1) in Hilbert Spaces}}},\ \href@noop {} {Ph.D. thesis},\ \bibinfo
  {school} {{University of Waterloo}} (\bibinfo {year} {{2011}}),\ \bibinfo
  {note} {retrieved from
  \url{https://uwspace.uwaterloo.ca/bitstream/handle/10012/6311/Hao_Yufang.pdf}}\BibitemShut
  {NoStop}%
\bibitem [{\citenamefont {Aleman}, \citenamefont {Martin},\ and\ \citenamefont
  {Ross}(2013)}]{AMR}%
  \BibitemOpen
  \bibfield  {author} {\bibinfo {author} {\bibfnamefont {A.}~\bibnamefont
  {Aleman}}, \bibinfo {author} {\bibfnamefont {R.}~\bibnamefont {Martin}}, \
  and\ \bibinfo {author} {\bibfnamefont {W.}~\bibnamefont {Ross}},\ }\bibfield
  {title} {\enquote {\bibinfo {title} {On a theorem of {L}ivsic},}\ }\href@noop
  {} {\bibfield  {journal} {\bibinfo  {journal} {J. Funct. Anal.}\ }\textbf
  {\bibinfo {volume} {264}},\ \bibinfo {pages} {999--1048} (\bibinfo {year}
  {2013})}\BibitemShut {NoStop}%
\bibitem [{\citenamefont {Clark}(1972)}]{Clark1972}%
  \BibitemOpen
  \bibfield  {author} {\bibinfo {author} {\bibfnamefont {D.}~\bibnamefont
  {Clark}},\ }\bibfield  {title} {\enquote {\bibinfo {title} {One dimensional
  perturbations of restricted shifts},}\ }\href@noop {} {\bibfield  {journal}
  {\bibinfo  {journal} {J. Anal. Math.}\ }\textbf {\bibinfo {volume} {25}},\
  \bibinfo {pages} {169--191} (\bibinfo {year} {1972})}\BibitemShut {NoStop}%
\bibitem [{\citenamefont {Krein}(1944{\natexlab{a}})}]{Krein1944}%
  \BibitemOpen
  \bibfield  {author} {\bibinfo {author} {\bibfnamefont {M.}~\bibnamefont
  {Krein}},\ }\bibfield  {title} {\enquote {\bibinfo {title} {On {H}ermitian
  operators with deficiency indices one},}\ }in\ \href@noop {} {\emph {\bibinfo
  {booktitle} {Dokl. Akad. Nauk SSSR}}},\ Vol.~\bibinfo {volume} {43}\
  (\bibinfo {year} {1944})\ pp.\ \bibinfo {pages} {339--342}\BibitemShut
  {NoStop}%
\bibitem [{\citenamefont {Krein}(1944{\natexlab{b}})}]{Krein1944one}%
  \BibitemOpen
  \bibfield  {author} {\bibinfo {author} {\bibfnamefont {M.}~\bibnamefont
  {Krein}},\ }\bibfield  {title} {\enquote {\bibinfo {title} {On one remarkable
  class of {H}ermitian operators},}\ }in\ \href@noop {} {\emph {\bibinfo
  {booktitle} {Dokl. Akad. Nauk SSSR}}},\ Vol.~\bibinfo {volume} {44}\
  (\bibinfo {year} {1944})\ pp.\ \bibinfo {pages} {191--195}\BibitemShut
  {NoStop}%
\bibitem [{\citenamefont {Silva}\ and\ \citenamefont
  {Toloza}(2007)}]{Silva2007}%
  \BibitemOpen
  \bibfield  {author} {\bibinfo {author} {\bibfnamefont {L.}~\bibnamefont
  {Silva}}\ and\ \bibinfo {author} {\bibfnamefont {J.}~\bibnamefont {Toloza}},\
  }\bibfield  {title} {\enquote {\bibinfo {title} {Applications of {K}rein's
  theory of regular symmetric operators to sampling theory},}\ }\href@noop {}
  {\bibfield  {journal} {\bibinfo  {journal} {J. Phys. A}\ }\textbf {\bibinfo
  {volume} {40}},\ \bibinfo {pages} {9413} (\bibinfo {year}
  {2007})}\BibitemShut {NoStop}%
\bibitem [{\citenamefont {Martin}(2011)}]{Martin-dB}%
  \BibitemOpen
  \bibfield  {author} {\bibinfo {author} {\bibfnamefont {R.}~\bibnamefont
  {Martin}},\ }\bibfield  {title} {\enquote {\bibinfo {title} {Representation
  of symmetric operators with deficiency indices $(1,1)$ in de
  \uppercase{B}ranges space},}\ }\href@noop {} {\bibfield  {journal} {\bibinfo
  {journal} {Complex Anal. Oper. Theory}\ }\textbf {\bibinfo {volume} {5}},\
  \bibinfo {pages} {545--577} (\bibinfo {year} {2011})}\BibitemShut {NoStop}%
\bibitem [{\citenamefont {deBranges}(1968)}]{dB}%
  \BibitemOpen
  \bibfield  {author} {\bibinfo {author} {\bibfnamefont {L.}~\bibnamefont
  {deBranges}},\ }\href@noop {} {\emph {\bibinfo {title} {Hilbert spaces of
  entire functions}}}\ (\bibinfo  {publisher} {Prentice Hall},\ \bibinfo {year}
  {1968})\BibitemShut {NoStop}%
\bibitem [{\citenamefont {Gorbachuk}\ and\ \citenamefont
  {Gorbachuk}(2012)}]{GG}%
  \BibitemOpen
  \bibfield  {author} {\bibinfo {author} {\bibfnamefont {M.}~\bibnamefont
  {Gorbachuk}}\ and\ \bibinfo {author} {\bibfnamefont {V.}~\bibnamefont
  {Gorbachuk}},\ }\href@noop {} {\emph {\bibinfo {title} {M.G. {K}rein's
  lectures on entire operators}}},\ Vol.~\bibinfo {volume} {97}\ (\bibinfo
  {publisher} {Birkh{\"a}user},\ \bibinfo {year} {2012})\BibitemShut {NoStop}%
\bibitem [{\citenamefont {Garcia}, \citenamefont {Mashreghi},\ and\
  \citenamefont {Ross}(2016)}]{GR-model}%
  \BibitemOpen
  \bibfield  {author} {\bibinfo {author} {\bibfnamefont {S.}~\bibnamefont
  {Garcia}}, \bibinfo {author} {\bibfnamefont {J.}~\bibnamefont {Mashreghi}}, \
  and\ \bibinfo {author} {\bibfnamefont {W.}~\bibnamefont {Ross}},\ }\href@noop
  {} {\emph {\bibinfo {title} {Introduction to model spaces and their
  operators}}},\ Vol.\ \bibinfo {volume} {148}\ (\bibinfo  {publisher}
  {Cambridge University Press},\ \bibinfo {year} {2016})\BibitemShut {NoStop}%
\bibitem [{\citenamefont {Hoffman}(2007)}]{Hoff}%
  \BibitemOpen
  \bibfield  {author} {\bibinfo {author} {\bibfnamefont {K.}~\bibnamefont
  {Hoffman}},\ }\href@noop {} {\emph {\bibinfo {title} {Banach spaces of
  analytic functions}}}\ (\bibinfo  {publisher} {Courier Corporation},\
  \bibinfo {year} {2007})\BibitemShut {NoStop}%
\bibitem [{\citenamefont {Granville}\ and\ \citenamefont
  {Martin}(2006)}]{granville2006prime}%
  \BibitemOpen
  \bibfield  {author} {\bibinfo {author} {\bibfnamefont {A.}~\bibnamefont
  {Granville}}\ and\ \bibinfo {author} {\bibfnamefont {G.}~\bibnamefont
  {Martin}},\ }\bibfield  {title} {\enquote {\bibinfo {title} {Prime number
  races},}\ }\href@noop {} {\bibfield  {journal} {\bibinfo  {journal} {The
  American Mathematical Monthly}\ }\textbf {\bibinfo {volume} {113}},\ \bibinfo
  {pages} {1--33} (\bibinfo {year} {2006})}\BibitemShut {NoStop}%
\bibitem [{\citenamefont {Oliver}\ and\ \citenamefont
  {Soundararajan}(2016)}]{oliver2016unexpected}%
  \BibitemOpen
  \bibfield  {author} {\bibinfo {author} {\bibfnamefont {R.~J.~L.}\
  \bibnamefont {Oliver}}\ and\ \bibinfo {author} {\bibfnamefont
  {K.}~\bibnamefont {Soundararajan}},\ }\bibfield  {title} {\enquote {\bibinfo
  {title} {Unexpected biases in the distribution of consecutive primes},}\
  }\href@noop {} {\bibfield  {journal} {\bibinfo  {journal} {Proceedings of the
  National Academy of Sciences}\ }\textbf {\bibinfo {volume} {113}},\ \bibinfo
  {pages} {E4446--E4454} (\bibinfo {year} {2016})}\BibitemShut {NoStop}%
\end{thebibliography}
\end{document}